\documentclass[12pt]{amsart}

 \newtheorem{thm}{Theorem}[section]
 \newtheorem{cor}[thm]{Corollary}
 \newtheorem{lem}[thm]{Lemma}
 \newtheorem{prop}[thm]{Proposition}
 \theoremstyle{definition}
 
 \theoremstyle{remark}
 \newtheorem{rem}[thm]{Remark}
 \theoremstyle{definition}
 \newtheorem{ex}[thm]{Example}

 \newcommand{\CC}{\mathbb{C}}

 \newcommand{\PP}{\mathbb{P}}

\usepackage{xypic} \usepackage{amsthm}

\begin{document}

\title{Complete intersection points on general surfaces in $\PP^3$}


\author[E. Carlini]{Enrico Carlini}
\address[E. Carlini]{Dipartimento di Matematica, Politecnico di Torino, Torino, Italia}
\email{enrico.carlini@polito.it}

\author[L. Chiantini]{Luca Chiantini}
\address[L. Chiantini]{Dipartimento di Scienze Matematiche e Informatiche, Universit\`{a} di Siena, Siena, Italia}
\email{chiantini@unisi.it}

\author[A.V. Geramita]{Anthony V. Geramita}
\address[A.V. Geramita]{Department of Mathematics and Statistics, Queen's University, Kingston, Ontario, Canada, K7L 3N6 and Dipartimento di Matematica, Universit\`{a} di Genova, Genova, Italia}
\email{Anthony.Geramita@gmail.com \\ geramita@dima.unige.it  }


\date{}


\begin{abstract}
In this paper we consider the existence of complete intersection
points of type $(a,b,c)$, on the generic degree $d$ surface of
$\PP^3$. For any choice of $a, \ b,\ c$ we resolve the existence question asymptotically, i.e. for all $d \gg 0$.  For small values of $a,\ b,\ c$ we resolve the existence problem completely.
\end{abstract}

\maketitle

\section{Introduction}

A recurrent theme in classical projective geometry is the study of
special subvarieties of some given family of varieties, e.g. how
many isolated singular points can a surface of degree $d$ in
$\PP^3$ have? when is it true that the members of a certain family
of varieties contain rational curves? contain a linear space of
some positive dimension?  Other examples of similar questions can be
easily provided by the reader.

The study of the special case of complete intersection
subvarieties of hypersurfaces in $\PP^n$ has been the subject of a great deal of research. It
was known to Severi \cite{Severi} that for $n \geq 4$ the only
complete intersections, of codimension one, on a general hypersurface are obtained by
intersecting that hypersurface with another.

This observation was extended to $\PP^3$ by Noether (and
Lefschetz)  \cite{Lefschetz, GH85} for general hypersurfaces of
degree $\geq 4$.  These ideas were further generalized by
Grothendieck \cite{SGA2}.

In \cite{CaChGe}, we proposed a new approach to the problem of
studying complete intersection subvarieties of hypersurfaces. This
approach used a mix of projective geometry and commutative algebra
and is more elementary and direct than, for example, the
approach of Grothendieck.  With our approach we were able to give a complete
description of the situation for complete intersections of
codimension $r$ in $\PP^n$ which lie on a general hypersurface of
degree $d$ whenever $2r \leq n+2$. The main result of
\cite{CaChGe} is the following:

\begin{thm}\label{nostro}
Let $X\subset\PP^n$ be a generic degree $d$ hypersurface, with
$n,d>1$. Then $X$ contains a complete intersection of type
$(a_1,\ldots,a_r)$, with $2r\leq n+2$, and the $a_i$ all less than
$d$, in the following (and only in the following) instances:
\medskip

\begin{itemize}
\item $n=2$: then $r=2$, $d$ arbitrary and $a_1$ and $a_2$ can
assume any value less than $d$;

\item $n=3$, $r=2 $:  for $d \leq 3$ we have that $a_1$ and $a_2$
can assume any value less than $d$;

\item $n=4$, $r = 3$: for $d \leq 5$ we have that $a_1, a_2$ and
$a_3$ can assume any value less than $d$;

\item $n=6,r = 4$ or $n=8,r = 5$: for $d \leq 3$ we have that
$a_1, \ldots ,a_r$ can assume any value less than $d$;

\item $n=5,7$ or $n>8$, $2r=n+1$ or $2r = n+2$:  we have only
linear spaces on quadrics, i.e. $d = 2$ and $a_1 = \ldots =a_r =
1$.

\end{itemize}

\end{thm}

In this paper, we are interested in the first case not covered by
Theorem \ref{nostro}. Namely, the case  $n=3,r=3$, i.e.
complete intersection points on surfaces of $\PP^3$. Although this
a very natural question, we are not aware of any reference to the
subject in the literature. Using the methods of \cite{CaChGe} we prove the
following:

\begin{thm}\label{mainTHM} For non-negative integers $a,b,c,d$, such that $a\leq b\leq c<
d$ we have the following:

\begin{itemize}

\item if $a\leq 4$, then the generic degree $d$ surface of $\PP^3$ contains a
$CI(a,b,c)$;

\item if $a=5,b\leq 11$, then the generic degree $d$
surface of $\PP^3$ contains a $CI(5,b,c)$; if $a=5,b=12$ and
$c=12$ then the generic degree $d$ surface of $\PP^3$ contains a
$CI(5,12,12)$; if $a=5,b=12$ and $c\geq 13$ then the generic
degree $d\geq 2c+15$ surface does not contain a $CI(5,12,c)$; if
$a=5$ and $b\geq 13$, then the generic degree $d\geq b+c+2$
surface does not contain a $CI(5,b,c)$;

\item if $a=6,b\leq 7$, then the generic degree $d$
surface of $\PP^3$ contains a $CI(6,b,c)$; if $a=6,b=8$ and
$c=8,9$ then the generic degree $d$ surface of $\PP^3$ contains a
$CI(6,8,c)$; if $a=6,\ b=8$ and $c\geq 10$ then the generic degree
$d\geq 2c+12$ surface does not contain a $CI(6,8,c)$; if $a=6$ and
$b\geq 9$, then the generic degree $d\geq b+c+3$ surface does not
contain a $CI(6,b,c)$;

\item if $a\geq 7$, then the generic degree $d\geq a+b+c-3$ surface of $\PP^3$
does not contain a $CI(a,b,c)$.
\end{itemize}
\end{thm}

Notice that Theorem \ref{mainTHM} gives a complete asymptotic solution to the existence problem for $CI(a,b,c)$ on a general surface of degree $d$ in $\PP^3$.  More precisely,

\begin{cor} Let $a \leq b \leq c <d$ be integers.  Then for $d\gg 0$ the generic degree $d$ surface contains a $CI(a,b,c)$ if:
\begin{itemize}
  \item $a \leq 4$;
  \item $a = 5,\ b\leq 11$;
  \item $a = 5,\ b=12,\ c=12$;
  \item $a=6,\ b\leq 7$;
  \item $a=6,\ b=8,\ c=8,9$.
 \end{itemize}
and does not contain a $CI(a,b,c)$ in all other cases.
\end{cor}

\begin{rem}
We also notice that the kind of asymptotic problem we solved above can only be considered for points.  More precisely, if we choose a family $\mathcal{F}$ of subschemes of $\PP^n$ we can ask the following: is it true that for $d\gg 0$ the generic degree $d$ hypersurface of $\PP^n$ contains an object of the family $\mathcal{F}$?

Using a standard incidence correspondence argument, it is easy to see that a positive answer can be given only if
$$
\dim {\mathcal{F}} + 1 - h_{\mathcal{F}}(d) \geq 0,
$$
where $h_{\mathcal{F}}(d)$ is the Hilbert polynomial of the objects in $\mathcal{F}$.  Clearly this can be the case only if $h_{\mathcal{F}}(d)$ is bounded and hence constant.  This implies that $\mathcal{F}$ is a family parameterizing $0$-dimensional schemes.
\end{rem}

The paper is structured as follows: in Section \ref{questionSEC},
we formalize the question we want to study and we treat the first
simple instances; in Section \ref{technicalFACTSsec}, we recall
the results we need from \cite{CaChGe}; in Sections
\ref{aleq4SEC}, \ref{ageq5nonexistSEC} and \ref{ageq5existSEC} we
apply our method to produce the intermediate results necessary to
prove Theorem \ref{mainTHM}. Finally, in Section \ref{finalSEC},
we prove Theorem \ref{mainTHM} and we state a conjecture for the
expected behavior in the cases which still remain open.

In the proof of Theorem \ref{ageq5THM} we used the
computer algebra system CoCoA \cite{cocoa} for which we thank the
developers of the software.

The first author wishes to thank Queen's University for its kind
hospitality during the writing of this paper, the research group
GNSAGA of INDAM and the special fund ``Fondo giovani ricercatori"
of the Politecnico di Torino for financial support.  The first and third author also wish to thank NSERC (Canada) for its financial support during the writing of this paper.

\section{The question}\label{questionSEC}
In this paper we study complete intersection points in projective
three space. We say that $\mathbb{X}\subset\PP^3$ is a complete
intersection 0-dimensional scheme if its ideal
$I_\mathbb{X}=(F,G,H)$ where the forms $F,G$ and $H$ are a regular
sequence in the ring $R=\mathbb{C}[x_0,\ldots,x_3]$. Moreover, if
$\deg F=a,\deg G=b$ and $\deg H=c$ we say that $\mathbb{X}$ is a
{\it complete intersection of type} $(a,b,c)$. We will always
assume $a\leq b\leq  c$ and we will write $CI(a,b,c)$
to describe a complete intersection of type $(a,b,c)$.

Our basic question is: {\em for which integers $a,b,c$ and $d$
does the general degree $d$ surface of $\PP^3$ contain a $CI(a,b,c)$}?

There are cases where the answer is straightforward. If $d=c$, the
answer is clearly affirmative as we are cutting a complete
intersection curve of type $(a,b)$ with a surface of degree $d$
(similarly for $d=a$ or $d=b$). If $d<a$, the answer is negative
as no form of degree less than $a$ belongs to the ideal of a
$CI(a,b,c)$, and similarly for $a<d<b$ as a generic form is
irreducible. If $b<d<c$, then we are really looking for a complete
intersection of type $(a,b)$ on the generic degree $d$ surface,
and this is dealt with in Theorem \ref{nostro}. Hence, it is
enough to focus on the following refinement of our question: {\em
for which integers $a,b,c$ and $d$, $a\leq b\leq c<d$, does the
general degree $d$ surface of $\PP^3$ contain a $CI(a,b,c)$}?

\section{Technical facts}\label{technicalFACTSsec}

We will treat this question using the method introduced in
\cite{CaChGe}. Our method proceeds as follows: translate the problem of finding a $CI(a,b,c)$ on a general surface of degree $d$, say $M = 0$, as the problem of writing $M$ as
$$
M = F F^\prime + G G^\prime + H H^\prime
$$
where $F,\ G,\ H$ and $F^\prime,\ G^\prime, \ H^\prime$ are forms of degree $a,\ b,\ c$ and $d-a,\ d-b,\ d-c$ respectively.  As $M$ is generic, this decomposition problem is actually a problem about joins of varieties of splitting forms.  Then we use Terracini's lemma to translate the computation of the dimension of the join, into a Hilbert function computation.  Namely, as first observed in \cite{Mamma}, the tangent space to the variety of splitting forms at the point $[FF^\prime]$ corresponds to the degree $d$ homogeneous piece of the ideal $ (F, F^\prime )$.  Thus, the tangent space at $M$ to the join corresponds to the degree $d$ homogeneous piece of the ideal spanned by $F,F^\prime, G, G^\prime, H, H^\prime$.  For more details we refer the reader to \cite{CaChGe}.

In particular we will need the following (see \cite[Lemma 4.3]{CaChGe}):

\begin{lem}\label{equivalence}
For given integers $a,b,c$ and $d$, such that $a\leq b\leq c<
d$, the following are equivalent facts:
\begin{enumerate}
\item  The general degree $d$ surface of $\PP^3$ contains a
$CI(a,b,c)$; \item\label{algebraicQUESTION} For a generic choice
of forms $F,G,H,H',G',F'\in R$ of degrees $a,b,c,d-c,d-b,d-a$ one
has that
    \[H\left({R\over (F,G,H,H',G',F')},d\right)=0\]
where $H(\cdot,d)$ denotes the Hilbert function in degree $d$.
\end{enumerate}
\end{lem}

Using Lemma \ref{equivalence} we translate our geometric question
into a purely algebraic one. In particular, we can take advantage of
results about the Lefschetz property \cite{Stanley,Anick} to deal
with our question.

As $F,G,H$ and $H'$ are a
regular sequence in $R$ we have a good understanding of the ring
\[W={R\over (F,G,H,H')}\]
and we will use this to study the Hilbert function of the ring
\[{R\over(F,G,H,H',G',F')}\simeq {W\over ([F'],[G'])},\]
where $[\cdot]$ denotes the class in $W$.

Via the Koszul complex we compute the minimal free resolution of $W$:
\begin{equation}\label{Wres}
0\leftarrow W\leftarrow M_0 \leftarrow M_1\leftarrow M_2\leftarrow
M_3\leftarrow M_4\leftarrow 0
\end{equation}
where
\[M_0=R,\]
\[M_1=R(-a)\oplus R(-b)\oplus R(-c)\oplus R(-d+c),\]
\[M_2=R(-a-b)\oplus R(-a-c)\oplus R(-a-d+c)\oplus R(-b-c)\oplus R(-b-d+c)\oplus R(-d),\]
\[M_3=R(-a-b-c)\oplus R(-a-b-d+c)\oplus R(-a-d)\oplus R(-b-d)\]
\[M_4=R(-a-b-d).\]

We also notice that (see  \cite[Lemma 4.1 and Remark 4.2]{CaChGe}):

\begin{lem}\label{switchLEM}
The following are equivalent:
\begin{itemize}
\item for integers $a\leq b\leq c<d$ a $CI(a,b,c)$ exists on the generic degree $d$ surface of $\PP^3$;
\item for integers $a'\leq b'\leq c'\leq {d}$ a $CI(a',b',c')$ exists on the generic degree $d$ surface of $\PP^3$, where $a=a'$ or $a+a'=d$, and $b=b'$ or $b+b'=d$, and $c=c'$ or $c+c'=d$.
\end{itemize}
\end{lem}

\begin{rem}\label{dover2LEM}
Using Lemma \ref{switchLEM} we can study our question for integers $a\leq b\leq c<d/2$ and produce a complete answer for the general case. In fact, either $a\leq d/2$ or $a'\leq d/2$.
\end{rem}


\section{The $a\leq 4$ case}\label{aleq4SEC}
Here we use Stanley's result \cite{Stanley} showing that the
quotient of $R=\mathbb{C}[x_0,\ldots,x_3]$ by four generic forms
has the Strong Lefschetz Property. More precisely, given generic
forms $F,G,H,F',G'\in R$, of degrees $a,b,c,d-c,d-b$, we consider
$W={R/ (F,G,H,H')}$. Then the multiplication by the class of $G'$
has maximal rank. Hence the sequence
\[W(-d+b){\rightarrow} W\rightarrow {W\over ( [G'])}\rightarrow 0\]
produces $H(W/([G'])),d)=\max\{H(W,d)-H(W,b),0\}$.

\begin{prop}\label{asymptoticamenodi5}
For any choice of $a,b,c$ and $d$ positive integers such that
$a\leq 4\leq b\leq c$ and $d\geq a+b+c-3$, the general degree $d$ surface in $\PP^3$
contains a $CI(a,b,c)$.

\end{prop}
\begin{proof} Using Lemma \ref{dover2LEM} it is enough to consider the case when $a\leq b\leq c\leq {d\over 2}$. Using Proposition \ref{equivalence} part \eqref{algebraicQUESTION} we only have to show that
\[H(W/([G']),d)=\max\{H(W,d)-H(W,b),0\}=0.\]
By the resolution of $W$ given in (\ref{Wres}) we immediately get:
\begin{itemize}
\item if $b<c$, then
\begin{eqnarray*}
H(W,b) & = & {b+3\choose 3}-{b-a+3\choose 3}-1 \\
       & = & 1/6 a^{3}-1/2a^{2}b + 1/2ab^{2} - a^{2} + 2ab + 11/6a-1;
\end{eqnarray*}
\item if $b=c$, then
\begin{eqnarray*}
H(W,b) & = & {b+3\choose 3}-{b-a+3\choose 3}-2 \\
       & = & 1/6 a^{3}-1/2a^{2}b + 1/2ab^{2} - a^{2} + 2ab + 11/6a-2.
\end{eqnarray*}
\end{itemize}

Led by the resolution of $W$, we also consider the following polynomial

\[h(W,d)=\]
\[={d+3\choose 3}-\left[{d-a+3\choose 3}+{d-b+3\choose
3}+{d-c+3\choose 3}+\right. \left. {c+3\choose 3}\right]+\]
\[+{d-a-b+3\choose
3}+{d-a-c+3\choose 3}+{c-a+3\choose 3}+{d-b-c+3\choose 3}+\]

\[+{c-b+3\choose3}+1-\left[{d-a-b-c+3\choose 3}\right. \left. +{c-a-b+3\choose 3}\right],\]
where ${x \choose 3}$ is the polynomial ${1\over 6}x(x-1)(x-2)$. Making the computation we get
\[h(W,d)=1/2 a^{2}b + 1/2ab^{2}-2ab + 1.\]

Notice that, for given $a,b$ and $c$ such that $c-a-b\geq -3$ and
$d\geq a+b+c-3$, the evaluation of $h(W,d)$ coincides with the
Hilbert function of $W$ in degree $d$, i.e. $h(W,d)=H(W,d)$.
Moreover, the inequalities
\[a\leq 4 \mbox{ and } c-a-b\leq -4\]
only hold when $a=4$ and $b=c$ (recall that $a\leq b\leq c$) and in this case
\[H(W,d)=h(W,d)+{c-a-b+3\choose 3}=h(W,d)-1.\]

Finally we compute $H(W,d)-H(W,b)$ distinguishing two cases.

The $a< 4$ or $b<c$ case. If $b<c$ we use the value of $H(W,b)$
and the equality $H(W,d)=h(W,d)$  previously determined to  get
\[H(W,d)-H(W,b)=-1/6 a^{3} + a^{2}b + a^{2}-4ab-11/6a + 2.\]
This polynomial is linear in $b$ and it does not involve $d$ and
it is easy to see that for $a\leq 4$
\[H(W,d)-H(W,b)\leq 0.\]
When $a < 4$ and $b=c$ a completely analogous argument can be applied.

The $a=4$ and  $b=c$ case. Mutatis mutandis, we compute again and we get
\[H(W,d)-H(W,b)=-1/6 a^{3} + a^{2}b + a^{2}-4ab-11/6a + 2,\]
hence the same polynomial of the previous case and this finishes the proof.
\end{proof}

Proposition \ref{asymptoticamenodi5} gives an asymptotic result
yielding that, when one of the degree of the $CI$ is at most 4, then
for $d$ big enough a complete intersection of the given type
exists on a generic surface of degree $d$. With a slightly more
careful analysis this can be improved and the condition on $d$ can
be dropped.

\begin{thm}\label{aleq4THM}
Let $a,b,c$ and $d$ be integers such that $a\leq  b\leq c<d$.
If $a\leq 4$, then a $CI(a,b,c)$ exists on the generic degree $d$
surface in $\PP^3$.
\end{thm}
\begin{proof}
Using Proposition \ref{asymptoticamenodi5} we have only to check
values of $d$ in the range $c<d\leq a+b+c-4$. The idea is to
use Lemma \ref{switchLEM} to reduce the degree of the complete
intersection not changing $d$ so that Proposition
\ref{asymptoticamenodi5} can be applied.

For $a=2$, we consider the existence of a $CI(2,b,c)$ on the
generic degree $d$ surface for $c<d\leq b+c-2$. Such a
complete intersection exists if the same happen for a
$CI(2,d-c,d-b)$. But, by Proposition \ref{asymptoticamenodi5},
this is the case as soon as
\[d\geq 2+(d-c)+(d-b)-3\]
and this equivalent to $d\leq b+c-1$ which is actually the case. Similarly for $a=3$.

The case $a=4$ is treated in analogy with the previous ones,
except for $d=b+c$. In this situation, applying Lemma
\ref{switchLEM}, we have to study $CI(4,b,b)$'s on a generic
surface of degree $d\geq b$. Repeating the same argument above we
have to treat values of $d$ in the range $b\leq d\leq 2b$.
Proposition \ref{existdbigPROP} gives the existence for all $d$,
but for $d=2b$. By Proposition $\ref{equivalence}$ we have to
consider the coordinate ring $W$ of a complete intersection of
type $(b,b,b,b)$ and its Hilbert function $H(\cdot)$. Using the
fact that multiplication by one form  has maximal rank in $W$, we
have only to compare $H(2b)$ and $H(2b-4)$, but these values are
the same being $W$ a Gorenstein ring with socle degree $4b-4$, and
this finishes the proof.
\end{proof}

\section{The case $a> 4$: non-existence results}\label{ageq5nonexistSEC}

In this section we will prove asymptotic non-existence results
when $a> 4$. For non-negative integers $a,b,c$ and $d$, $a\leq
b\leq c<d$, we consider generic forms $F,G,H,H',G',F'\in R$ of
degrees $a,b,c,d-c,d-b$ and $d-a$. Consider the ring $W={R/
(F,G,H,H')}$ and notice that, by a straightforward dimensional
argument, if
\[H(W,a)+H(W,b)-H(W,d)<0\]
then
\[H\left({W\over ([G'],[F'])},d\right)\neq 0.\]
Hence, by Lemma \ref{equivalence} \eqref{algebraicQUESTION},  if
$H(W,d)-H(W,a)-H(W,b)<0$, then the generic degree $d$ surface of
$\PP^3$ does not contain a
$CI(a,b,c)$. Using this idea we prove the following:

\begin{thm}\label{nonEXISTa5}
Let $a\leq b\leq c$ and $d$  be non-negative integers such that
\[a=5 \mbox{ and } b\geq 13\]
or
\[a=6 \mbox{ and } b\geq 9\]
or
\[a \geq 7.\]

Then, for $d\geq a+b+c-3$ the generic degree $d$ surface of
$\PP^3$ does not contain a
$CI(a,b,c)$.
\end{thm}

In order to prove this theorem we need the following technical
fact:

\begin{lem}\label{inequalityLEMMA}
Let $a,b,c$ be non-negative integers, such that $4<a$ .
Assume that, for integers $c_0$ and $d$ such that
\[c_0\geq b \mbox{ and } d>a+b+c_0-4,\]
one has the Hilbert function inequality
\[H(W,a)+H(W,b)-H(W,d)<0,\]
where $W$ is the ring
\[
W={R\over (F,G,H,H')}
\]
and the forms $F,G,H$ and $H'$ are generic and have degrees
$a,b,c_0$ and $d-c_0$.

Then, if $A,\ B,\ C$ and $D$ are forms of degrees $a,\ b,\ c\geq c_0$ and $d-c$ and $$W^\prime = {R\over (A,B,C,D)}$$ then the following inequality holds:
\[H(W^\prime,a)+H(W^\prime,b)-H(W^\prime,d)<0,\]
for $d>a+b+c-4$.
\end{lem}
\begin{proof}
The key observation is that
\[H(W,d)=H(W,a+b-4).\]
In fact, being $W$ a Gorenstein ring,  its Hilbert function is symmetric and $H(W,x)=H(W,y)$ if $x+y=d+a+b-4$.
Then, we compute
\[H(W,a)+H(W,b)-H(W,a+b-4)\]
using the formulae in the proof of Proposition
\ref{asymptoticamenodi5}, for which we need the assumption on $d$.
One sees that the final expression does not involve neither $c$ or
$d$ and the proof follows. For example, in the case $a<b$ one gets
$H(W,a)+H(W,b)-H(W,a+b-4)<0$ if and only if
\[b\geq \frac{{1\over 2}a^3-a^2+{11\over 2}a -4}{(a^2-4a)}.\]
\end{proof}

We can now prove Theorem \ref{nonEXISTa5}.

\begin{proof}[Proof of Theorem \ref{nonEXISTa5}]
We let $b=c$ and we show that for the required values of $a$ and
$b$ we have the inequality $H(W,a)+H(W,b)-H(W,d)<0$. Then we apply
Lemma \ref{inequalityLEMMA} to get the result when $c\geq b$.

Again, we notice that $H(W,d)=H(W,a+b-4)$.

We divide the proof in two cases depending on whether $a=b$ or
$a<b$.

\noindent{\bf Case $a<b$.}

Using the resolution of the ring $W$, the
inequality
\[H(W,a)+H(W,b)-H(W,a+b-4)\geq 0\]
is readily seen to be equivalent to
\[b\leq \frac{{1\over 2}a^3-a^2+{11\over 2}a -4}{(a^2-4a)}.\]
Recalling that $a<b$ we get
\[H(W,a)+H(W,b)-H(W,d)\geq 0\]
only if
\[-{1\over 2}a^3+3a^2+{11\over 2}a-4\geq 0\]
and this inequality holds if and only if
\[a=5 \mbox{ or } a=6.\]
Hence
\[H(W,a)+H(W,b)-H(W,d)\geq 0\]
implies
\[a=5,b\leq 12\]
or
\[a=6,b\leq 8.\]

{\bf Case $a=b$.}

 Computing we get
\[H(W,a)+H(W,b)-H(W,d)=-{2\over 3}a^3+4a^2+{11\over 3}a-3\geq 0\]
only if $a< 7$ and this finishes the proof.
\end{proof}

To prove some more non-existence results, we need the following:

\begin{prop}\label{nonexistEVENTUALLY}

Let $a\leq b\leq c< {d}$ and $d> 2c+b+a-3$. If no $CI(a,b,c)$
exists on the generic degree $d$ hypersurface, then it does not
exist on the generic hypersurface of degree $d'>d$ either.
\end{prop}
\begin{proof}
It is enough to treat the case $d'=d+1$.  For generic forms $F,\ G,\ H$ of degrees $a,\ b$ and $c$ let
$$
A = {\CC[x_0,\cdots,x_3]\over(F,G,H)}.
$$
By hypothesis, for the
generic choice of $F',G'$ and $H'$ of degrees $d-a,d-b$ and $d-c$
in $A$, we know that the degree $d$ part of
\[{A\over (F',G',H')}\]
is not zero. Now, consider elements $F'',G''$ and $H''$ of degrees
$d+1-a,d+1-b$ and $d+1-c$. Notice that

\[d-a\geq d-b\geq d-c >c+b+a-3\]
and recall that $A_i\simeq A_j$ as $\CC$ vector spaces if $i$ and $j$ are $>2c + b + a - 3$. Thus for a general linear form
$L$ we have
\[F''=LF^*, G''=LG^*\mbox{ and } H''=LH^*\]
and the forms $F^*,G^*$ and $H^*$ have degrees $d-a,d-b$ and
$d-c$. Hence we have a isomorphism
\[(F'',G'',H'')_{d+1}\simeq (F^*,G^*,H^*)_d\]
and this is enough to conclude that the degree $d+1$ part of
\[{A\over (F'',G'',H'')}\]
is not zero and the result follows.
\end{proof}

\begin{lem}\label{nonexistb8and12}
If $c\geq 13$, then the generic degree $d\geq 2c+15$ surface of
$\PP^3$ does not contain a $CI(5,12,c)$.

If $c\geq 10$, then the generic degree $d\geq 2c+12$ surface of
$\PP^3$ does not contain a $CI(6,8,c)$.
\end{lem}

\begin{proof}
We begin with the study of $CI(5,12,c)$. Let $c=13+x$,
$d=2c+a+b-2=41+2x$ and consider the ring
\[
W={R\over (F,G,H,H')}
\]
where the forms $F,G,H$ and $H'$ have degrees $5,12,13+x$ and
$28+x$. The generic degree $d$ surface does not contain a
$CI(5,12,c)$ if $H(W,5)+H(W,12)-H(W,d)<0$, where
$H(W,d) = H(W,41+2x)=H(W,14)$. Now we compute
\[H(5)={8\choose 3}-1=55,\]
\[H(12)={15\choose 3}-{10\choose 3}-1=334,\]
\[H(14)=\left\lbrace
\begin{array}{lr}
{17\choose 3}-{12\choose 3}-{5\choose 3}=450 & \mbox{ if } x>1 \\ \\
449 & \mbox{ if } x=1 \\ \\
446 & \mbox{ if } x=0
\end{array}
\right. .\] Hence, $H(W,5)+H(W,12)-H(W,d)<0$, and by Proposition
\ref{nonexistEVENTUALLY} we conclude that the generic degree $d'$
surface does not contain a $CI(5,12,c)$ for $d'\geq
d=41+2x=15+2c$.

The case of $CI(6,8,c)$ is solved by completely analogous computations.
\end{proof}

\section{The case $a> 4$: existence results}\label{ageq5existSEC}

Theorem \ref{nonEXISTa5} does not cover small values of $a$ and $b$.
In this Section we derive a result analogous to Theorem
\ref{aleq4THM} in these cases.

We begin with proving two technical facts.

\begin{prop}\label{existdbigPROP}
Let $a\leq b\leq c\leq {d}$ and $d\geq a+b+c-3$. If a $CI(a,b,c)$
exists on the generic degree $d$ surface, then it also exists on
the generic surface of degree $d'>d$.
\end{prop}
\begin{proof}
Let $d'=d+1$ and notice that it is enough to treat this case. The
hypothesis reads as follows: the degree $d$ part of the ring
\[{A\over (F',G',H')}\]
is zero for generic forms $F',G'$ and $H'$ of degrees $d-a,d-b$
and $d-c$ where $A=R/(F,G,H)$. If $L\in A$ is a generic linear
form, by \cite{Stanley}, we know that multiplication by $L$ is
an isomorphism in degree bigger than or equal to $a+b+c-4$. Hence,
the degree $d+1$ piece of
\[{A\over (LF',LG',LH')}\]
is zero and this is enough to complete the proof since, if three special
forms, namely $LF',LG'$ and $LH'$, have maximal span then the same property holds
for a generic choice.
\end{proof}

\begin{lem}\label{existdbigLEM}
Let  $a,b$ and $d$ be non-negative integers such that $4<a\leq b$
and $d=2a+2b-6$. If the generic degree $d$ surface in $\PP^3$
contains a $CI(a,b,a+b-3)$, then the generic degree $d'$ surface
contains a $CI(a,b,c)$ for any $d'\geq a+b+c-3$ and any $c\geq
a+b-3$.
\end{lem}
\begin{proof}
Notice that, by Proposition \ref{existdbigPROP}, the generic
degree $d+\epsilon$ surface in $\PP^3$ contains a $CI(a,b,a+b-3)$
for all $\epsilon\geq 0$. Hence, by Lemma \ref{switchLEM}, the
same holds for $CI(a,b,a+b-3+\epsilon)$ and surfaces of degree
$d+\epsilon$. Making $\epsilon$ vary and again applying
Proposition \ref{existdbigPROP} the result follows.
\end{proof}

\begin{thm}\label{ageq5THM}
Let  $a,b,c$ and $d$ be non-negative integers such that $a\leq
b\leq c<d$. If $a=5$ and $b\leq 11$, or $a=6$ and $b\leq 7$,
then a $CI(a,b,c)$ exists on the generic degree $d$
surface of $\PP^3$. If $a=5,b=12$ and $c=12$, or $a=6,b=8$ and $c=8,9$,
then a $CI(a,b,c)$ exists on the generic degree $d$
surface of $\PP^3$.
\end{thm}

\begin{proof}
To prove the thesis we combine all the previous results and
technical facts. Crucial ingredients are also some explicit
computations that we performed using the computer algebra system
CoCoA \cite{cocoa}.

To determine whether a $CI(a,b,c)$ exists on the generic surface of degree $d$ in $\PP^3$, we proceed as follows:
\begin{itemize}
\item if $c\leq a+b-3$, we make explicit computations for all
$d\leq a+b+c-3$; a positive answer for $d=a+b+c-3$ solves the
cases for bigger $d$'s by Proposition \ref{existdbigPROP};
\item if $c=a+b-3$ and $d\geq 2a+2b-6$, we verify each statement
with an explicit computation for $d=2a+2b-6$; if the answer is
positive we conclude the same for $c\geq a+b-3$ and $d\geq
a+b+c-3$ by Lemma \ref{existdbigLEM}.
\end{itemize}

We sketch this procedure for $a=6$, the case $a=5$ is completely
analogous but lengthier. We need to perform explicit computations in the following cases:
\begin{itemize}
\item $CI(6,6,c_1)$, for $c_1\leq 9$ and $d_1\leq 9+c_1$;
\item $CI(6,7,c_2)$, for $c_2\leq 10$ and $d_3\leq 10+c_2$;
\end{itemize}
The computations (see Example \ref{computationEX}) show that the  complete intersections exist on
the generic surfaces of the required degrees. Hence we conclude
that the generic surface of degree $d$ contains a $CI(6,b,c)$ for
all $b\leq 7$ and any $c,d$ such that $d> c$. We conclude the proof for $a=6$ by verifying existence in the cases: $CI(6,8,8)$ for $d=19$, and $CI(6,8,9)$ for $d=20$.

\end{proof}

\begin{ex}\label{computationEX}
We begin with verifying that the generic surface of degree $7\leq d\leq 15$ contains a
$CI(6,6,6)$. Using Proposition \ref{equivalence} it is enough to
show that the ring
\[S={\CC[x_0,\ldots,x_3]\over(F,G,H,H',G',F')}\]
is zero in degree $d$, where the forms $F,G,H,H',G'$ and $F'$ are
generic and have degrees $6,6,6,d-6,d-6$ and $d-6$. Hence, for each $d$, we choose
random forms with rational coefficients of the  required degrees.
Then we ask CoCoA \cite{cocoa} to compute the Hilbert function of
$S$ in degree $d$. Since for all $d$'s we get $H(S,d)=0$, we conclude (by semicontinuity)
that this is the case for a generic choice of forms of the
appropriate degrees. In particular, as $15=6+6+6-3$ and $H(S,15)=0$, Proposition \ref{existdbigPROP} yields that a $CI(6,6,6)$ exists on the generic degree $d\geq 15$ surface of $\PP^3$.

The same argument works in complete analogy for $c\leq 8$. For $c=9$ we make an explicit computation for $d=18$ and using Lemma \ref{existdbigLEM} we show existence of a $CI(6,6,c)$ on the generic degree $d$ surface for $c\geq 9$ and $d\geq c+9$. The cases for $c<d<c+9$ are solved using Lemma \ref{equivalence} and the results for $c\leq 8$ and $a\leq 4$.
\end{ex}

\section{Main theorem and final remarks}\label{finalSEC}

We can now prove our main theorem:
\begin{proof}[Proof of Theorem \ref{mainTHM}]
The existence part for the case $a\leq 4$ is Theorem
\ref{aleq4THM} while existence for the remaining cases is Theorem
\ref{ageq5THM}. The asymptotic non-existences are given by Lemma
\ref{nonexistb8and12} and Theorem \ref{nonEXISTa5}.
\end{proof}

Theorem \ref{mainTHM} produces a complete asymptotic answer to our
original question. We also get many existence and non existence
results for small value of $d$. However, there are still
infinitely many cases which we have not solved, e.g. $a=7$ any
$b,c$ and $d$ such that $7\leq b\leq c$ and $c+5\leq d\leq
a+b+c-4$.

We state a conjecture completing Theorem \ref{mainTHM}:
\begin{quote}{\it{ \textbf{Conjecture}:}
given non-negative integers $a,b,c$ and $d$ such that $ a\leq
b\leq c<d$, there exists a function $d(a,b,c)$, possibly assuming
the value $+\infty$, such that the generic degree $d$ surface in
$\PP^3$ contains a $CI(a,b,c)$ if and only if $d<d(a,b,c)$.
}\end{quote}

As support for this conjecture, notice that it fits with the
asymptotic statement and with the other results of Theorem
\ref{mainTHM}. For example, $d(a,b,c)=+\infty$ for $a\leq 4$ and
$d(a,b,c)<a+b+c-3$ for $7\leq a$.




\bibliographystyle{alpha}

\def\cprime{$'$}

\end{document}